\documentclass[12pt]{article}
\usepackage{amssymb}
\usepackage{amsmath}
\usepackage{theorem}
\usepackage{latexsym}
\theoremheaderfont{\bf} 
\theorembodyfont{\sl}
\newtheorem{theo}{Theorem}[section]
{\theorembodyfont{\rm} \newtheorem{defi}[theo]{Definition}}
{\theorembodyfont{\rm} }
{\theorembodyfont{\rm} \newtheorem{rem}[theo]{Remark}}
{\theorembodyfont{\rm} \newtheorem{notation}[theo]{Notation}}
{\theorembodyfont{\rm} }
\newtheorem{prop}[theo]{Proposition}
\newtheorem{cor}[theo]{Corollary}
\newtheorem{lemma}[theo]{Lemma}

\newenvironment{proof}{{\sc Proof:}}{\mbox{}\hfill$\Box$\par}
\newcommand{\eqr}[1]{~\mbox{$(${\rm \ref{#1}}$)$}}
\newcommand{\Section}[1]{\section{#1}\setcounter{equation}{0}}

\newcommand{\junk}[1]{}

\newcommand{\N}{{\mathbb N}}
\newcommand{\F}{{\mathbb F}}

\newcommand{\Z}{{\mathbb Z}}
\newcommand{\C}{{\mathcal C}}

\newcommand{\wt}{{\rm wt}}
\newcommand{\rank}{{\rm rank}\,}

\newcommand{\dfree}{\mbox{$d_{\mbox{\rm\tiny free}}$}}

\newcommand{\zwei}[2]{\left[ \begin{array}{c}
                   #1 \\ #2 \end{array} \right]}

\newcounter{abc}

\title{Convolutional Codes with Maximum Distance Profile\footnote
    {This work was supported in part by NSF grants DMS-00-72383 and 
    CCR-02-05310.}}
\date{July 14, 2003}
\author{
  Ryan Hutchinson\\
  {\small Department of Mathematics}\vspace{-2mm}\\
  {\small University of Notre Dame}\vspace{-2mm}\\
  {\small Notre Dame, IN 46556-5683, USA}\vspace{-2mm}\\
  {\small {\em e-mail:} rhutchin@nd.edu} \and
  Joachim Rosenthal\\
  {\small Department of Mathematics}\vspace{-2mm}\\
  {\small University of Notre Dame}\vspace{-2mm}\\
  {\small Notre Dame, IN 46556-5683, USA}\vspace{-2mm}\\
  {\small {\em e-mail:} Rosenthal.1@nd.edu} \and
  Roxana Smarandache\\
  {\small Department of Mathematics and Statistics}\vspace{-2mm}\\
  {\small San Diego State University}\vspace{-2mm}\\
  {\small 5500 Campanile Dr.}\vspace{-2mm}\\
  {\small San Diego, CA 92182-7720, USA}\vspace{-2mm}\\
  {\small {\em e-mail:} rsmarand@sciences.sdsu.edu }}

\begin{document}
\maketitle
\begin{abstract}
  Maximum distance profile codes are characterized by the
  property that two trajectories which start at the same state
  and proceed to a different state will have the maximum possible
  distance from each other relative to any other convolutional
  code of the same rate and degree.

  In this paper we use methods from systems theory to
  characterize maximum distance profile codes algebraically. The
  main result shows that maximum distance profile codes form a
  generic set inside the variety which parameterizes the set of
  convolutional codes of a fixed rate and a fixed degree.\bigskip

\end{abstract}

\clearpage
\Section{Introduction}

The concept of maximum distance profile codes was introduced
in~\cite{gl03r}. This concept is closely related to the concept
of optimum distance profile code widely studied in the
convolutional code area; see e.g.~\cite{jo89,jo99}.  It was shown
in~\cite{gl03r} that maximum distance profile codes exist when
the transmission rate is $(n-1)/n$, and it was conjectured that
such codes exist for every transmission rate. In systems
theoretic terms, this means that the existence of maximum
distance profile codes was established for multi-input,
one-output systems and that it was conjectured for general
multi-input, multi-output (MIMO) systems.

The main result of this paper shows the existence of maximum
distance profile codes for general MIMO systems and that these
codes are generic in the sense of algebraic geometry. The
techniques we are using to establish this result are based on
very classical results from linear systems theory. Thus, we will
first explain the problem in terms of linear systems theory.  For
this, we follow the description as it can be found
in~\cite{ro96a1,ro99a}.

Let $\F$ be a finite field.  Let $n,k$, and $\delta$ be positive
integers with $k<n$.  Consider the matrices $A \in \F^{\delta
  \times \delta}, \,\,\, B \in \F^{\delta \times k}, \,\,\, C \in
\F^{(n-k) \times \delta}$, and $D \in \F^{(n-k) \times k}$.  A
rate $k/n$ convolutional code $\C$ of degree $\delta$ can be
described by the linear system governed by the equations:
\begin{eqnarray}  \label{iso}
    x_{t+1} & = & Ax_t+Bu_t, \nonumber \\  y_t & = & Cx_t+Du_t, \\
   v_t & = &\binom{y_t}{u_t}, \,\,  x_0=0. \nonumber
\end{eqnarray}
We call $x_t\in\F^{\delta}$ the {\em state vector},
$u_t\in\F^{k}$ the {\em information vector}, $y_t\in\F^{n-k}$ the
{\em parity vector}, and $v_t\in\F^n$ the {\em code vector}, each
at time~$t$.  The set of all possible code vector sequences
$v_t\in\F^n$ is called the convolutional code generated by
$(A,B,C,D)$; its elements are called codewords. 
For simplicity, we will assume that $(A,B)$ forms a
controllable pair and $(A,C)$ forms an observable pair. We will
refer to such a code as an $(n,k,\delta)$-code.

By iterating the equations defining the system~(\ref{iso}), it
can be seen that a sequence $\{v_t = \binom{y_t}{u_t}\in\F^n
~\mid ~t=0,1,2,\ldots ,j\}$ represents the beginning of a
codeword 
if and only if the following matrix equation is
satisfied:
\begin{equation}                                        
\left(
\begin{array}{ccc|ccccc}
&  & &      D        &               &      &  &  \\
&  & &     CB        &        D      &      &  &  \\
&-I&\ \ &     CAB    &        CB     &\ddots&  &  \\
&  & &        \vdots &               &\ddots&\ddots &  \\
&  & &CA^{j -1}B&CA^{j -2}B&\cdots& CB & D
\end{array}
\right) \left(
\begin{array}{c}
y_0 \\
y_1 \\
\vdots  \\
y_j  \\   \hline
u_0\\
u_1 \\
\vdots  \\
u_j
\end{array}
\right) =0.  \label{nparity}
\end{equation}

For the purpose of error control coding, it is important that any 
two codewords are far apart with respect to a suitable metric.
The following definition is fundamental in coding theory:
\begin{defi}
  Let $x,y\in \F^n$ be vectors.  The Hamming distance Ham($x,y$)
  is defined to be the number of components in which $x$ and $y$
  differ.  The weight wt($x$) of $x$ is defined to be the number
  of nonzero components of $x$.
\end{defi}
Clearly, one has that Ham($x,y$) = wt($x - y$). When $\{v_t =
\binom{y_t}{u_t}\in\F^n~\mid ~t=0,1,2,\ldots\}$ is a codeword,
one defines its weight to be $\sum_t \mbox{wt(} v_t\mbox{)}$.

In this paper, we are concerned only with finite-weight
codewords.  These are defined as follows:
\begin{defi}
  A sequence $\{v_t = \binom{y_t}{u_t}\in\F^n ~\mid
  ~t=0,1,2,\ldots\}$ represents a $finite$-$weight\\~codeword$ if
  \begin{enumerate}
  \item Equation~(\ref{iso}) is satisfied for all $t\in\Z_+$,
    where $\Z_+$ denotes the set of positive integers;
  \item There is an integer $j$ such that $x_{j+1} = 0$ and $u_t
    = 0$ for $t\geq j + 1$.
  \end{enumerate}
\end{defi}

A well-studied concept in convolutional coding theory is that of
column distances~\cite{jo99}. We give a systems theoretic
definition.

\begin{defi} 
  The weight wt($v$) of a vector $v$ is the number of nonzero 
  components of $v$.
\end{defi}

\begin{defi}
  The $j$th column distance of the code $\C$ is defined as
  $$
  d_j:=\min\left\{\sum_{t=0}^j\wt(u_t)+
    \sum_{t=0}^j\wt(y_t)\right\} ,
  $$
  where the minimum is taken over all trajectories $(u_t,y_t)$
  of the system\eqr{iso} with initial vector $u_0\neq 0$.
\end{defi}

Clearly, one has that $d_0\leq d_1\leq d_2 \leq \ldots $, and
hence there exists an integer $r$ such that $d_r=d_{r+j}$ for
all $j\geq 0$. This largest possible column distance is of
central importance in coding theory:

\begin{defi}
  \begin{equation}
    \dfree:=\lim_{j\rightarrow \infty} d_j
  \end{equation}
  is called the {\em free distance} of the code $\C$.
\end{defi}

Codes with a large free distance and the largest possible column
distances are very desirable. The following two results give
estimates for these parameters. The first one was proven
in~\cite{gl03r}:
\begin{prop}                  \label{P-dcj.bound}
  For every $j\in\N_0$, we have
\[
d_j\leq(n-k)(j+1)+1.
\]
\end{prop}
The proof given in~\cite{gl03r} uses algebraic properties of the
parity check matrix. The following systems theoretic proof is
almost trivial.

\noindent\begin{proof}
  Take an input sequence $u_0,\ldots,u_j$, where $\wt(u_0)=1$ and
  $u_s=0$ for $s\geq 1$. Let $y_0,\ldots,y_j\in \F^{n-k}$ be the
  corresponding output sequence. Then,
  $$
  d_j\leq \sum_{t=0}^j\wt(u_t)+ \sum_{t=0}^j\wt(y_t)\leq
  (n-k)(j+1)+1.
  $$
\end{proof}

The following theorem gives an upper bound for the free distance.
\begin{theo}[\cite{ro99a1}]               \label{T-MDS}
  The free distance of an $(n,k,\delta)$-code satisfies
\begin{equation}                                \label{G-Singleton}
   \dfree\leq (n-k)\Big(\Big\lfloor\frac{\delta}{k}\Big\rfloor+1\Big)
   +\delta+1.
\end{equation}
\end{theo}
The bound on the right hand side is called the {\em generalized
  Singleton bound}.  With these preliminaries, we can now give
the following definitions:
\begin{defi} 
  Let $\C$ be an $(n,k,\delta)$-code with column distances $d_j$
  and free distance $\dfree$.
  \begin{enumerate}
  \item $\C$ is said to have a {\em maximum distance profile} if
    $$
    d_j=(n-k)(j+1)+1\mbox{ for }j=0,\ldots,
    L:=\Big\lfloor\frac{\delta}{k}
    \Big\rfloor+\Big\lfloor\frac{\delta}{n-k}\Big\rfloor.
    $$
  \item $\C$ is called an {\em MDS code} if $\dfree$ attains the
    generalized Singleton bound\eqr{G-Singleton}.
  \item $\C$ is called a {\em strongly MDS code} if
   \[
   d_M=(n-k)\Big(\Big\lfloor\frac{\delta}{k}
   \Big\rfloor+1\Big)+\delta+1 \text{ for }
   M=\Big\lfloor\frac{\delta}{k}\Big\rfloor+
   \Big\lceil\frac{\delta}{n-k}\Big\rceil.
    \]
  \end{enumerate}
\end{defi}
In~\cite{ro99a1,sm01a}, it was shown that, for any rate
$k/n$ and degree $\delta$, MDS codes form a generic set
in the variety parametrizing convolutional codes of rate $k/n$
and degree $\delta$. In~\cite{gl03r},
the existence of $(n,n - 1,\delta)$ strongly MDS codes was 
established. When $n-k$
divides $\delta$, we have $M=L$.  It follows that, in this
situation, a convolutional code has a maximum distance profile if
and only if it is strongly MDS.  In Theorem~\ref{main} of this
paper, 
we establish the
existence of maximum distance profile convolutional codes for all 
parameters $(n,k,\delta)$ over a suitably large base field $\F$.

\Section{Existence of Maximum Distance Profile Codes}

The set of
all $(n,k,\delta)$ convolutional codes has in a natural way the
structure of a quasi-projective variety. For this, note that the
set of 4-tuples $(A,B,C,D)$, with $A \in \F^{\delta \times
  \delta}, \,\,\, B \in \F^{\delta \times k}, \,\,\, C \in
\F^{(n-k) \times \delta}$, $D \in \F^{(n-k) \times k}$, $(A,B)$ a
controllable pair, and $(A,C)$ an observable pair, describes the
set of $(n-k)\times k$ proper transfer functions of McMillan
degree $\delta$.  Hazewinkel~\cite{ha77p} showed that this set is
not only a quasi-projective variety but even a quasi-affine
variety. We may also view this set as a Zariski open subset of
the projective variety described in~\cite{ra94}. In this section, 
we establish the existence of maximum distance
profile codes for all parameters $(n,k,\delta)$ for sufficiently
large fields.  Moreover, we show
that the set of maximum distance profile codes forms a {\em
  generic set} when viewed as a subset of the quasi-projective
variety of all $(n,k,\delta)$ convolutional codes.  More
precisely, we show that the set of maximum distance profile
codes is open and dense inside this quasi-projective variety.

The strategy for obtaining this result is as follows. In the
first step (Theorem~\ref{main1}), we exhibit a set of polynomial
equations whose zero set exactly describes the $(n,k,\delta)$
codes which do {\bf not} have the maximum distance profile
property. This shows that the codes possessing the maximum
distance profile property form a Zariski open subset.  In the
second step, we show that this Zariski open subset is nonempty as
soon as the field is sufficiently large. This part of the proof
invokes some classical results from partial realization theory.

The block Toeplitz matrix appearing in\eqr{nparity} is
of central importance in what follows.  Thus, we define:

\begin{equation}                          \label{Bl-To}
       \mathcal{T}_j  := 
  \left(
     \begin{array}{cccc}
                     F_0    &          &         &         \\
                     F_1    & F_0      &         &         \\
                     \vdots & \vdots   & \ddots  &         \\
                     F_j    & F_j-1    & \cdots  & F_0
     \end{array}
  \right)
:=
 \left(
  \begin{array}{ccccc}
  D           &                    &          &       &  \\
  CB          &        D           &          &       &  \\
  CAB         &        CB          &\ddots    &       &  \\
  \vdots      &                    &\ddots    &\ddots &  \\
  CA^{j -1}B  &        CA^{j -2}B  &\cdots    & CB    & D
  \end{array}
  \right).
\end{equation}

\begin{notation}     \label{B-minors}
  Let $i_1<\ldots<i_r\leq (j+1)(n-k)$ and $j_1<\ldots<j_r\leq
  (j+1)k$ be two sets of indices. We denote by
  $M^{i_1,\ldots,i_r}_{j_1,\ldots,j_r}$ the $r\times r$ minor
  obtained from~$\mathcal{T}_j$ by picking the rows with indices
  $i_1,\ldots,i_r$ and the columns with indices $j_1,\ldots,j_r$.
\end{notation}

It will turn out that an $(A,B,C,D)$ code has $j$th column
distance $d_{j}=(n-k)(j+1)+1$ if and only if all minors appearing
in\eqr{Bl-To} which are not trivially zero are nonzero. In order
to make this statement precise, we make the following definition.
In this definition, we think of the nonzero entries of the block
Toeplitz matrix $\mathcal{T}_j\in \F^{(j+1)(n-k)\times (j+1)k}$ as
indeterminates of the polynomial ring
$R:=\F[x_1,x_2,\ldots,x_{(j+1)(n-k)k}]$.  Specifically, if the
entry $(s,t)$ of the matrix $F_i$ is nonzero, we set it
equal to $x_{i(n - k)k + (s - 1)k + t}$; otherwise, we leave it zero.

\begin{defi} 
  A minor $M^{i_1,i_2,\ldots,i_r}_{j_1,j_2,\ldots,j_r}$ of
  $\mathcal{T}_j$ is called {\em trivially zero} if
  $M^{i_1,i_2,\ldots,i_r}_{j_1,j_2,\ldots,j_r}$ is zero when
  viewed as an element of the ring $R$ in the manner specified in 
  the preceding paragraph.
\end{defi}
The following Lemma gives an algebraic criterion for a minor to
be trivially zero.
\begin{lemma}                 \label{B-minors2}
  Let $\mathcal{T}_j\in \F^{(j+1)(n-k)\times (j+1)k}$ be a block
  Toeplitz matrix as defined above. Suppose that, for all
  integers $i$ with
  $0\leq i\leq j$, all entries of the matrix $F_i$ are nonzero.  
  Then, a minor
  $M^{i_1,i_2,\ldots,i_r}_{j_1,j_2,\ldots,j_r}$ of
  $\mathcal{T}_j$ is not trivially zero if and only if its
  indices $ i_1<\ldots<i_r\leq (j+1)(n-k)$ and $
  j_1<\ldots<j_r\leq (j+1)k$ satisfy
  \begin{equation}\label{eq-indices}
     j_t\leq \Big\lceil\frac{i_t}{n-k}\Big\rceil k\mbox{ for }
   t=1,\ldots ,r.
   \end{equation}
\end{lemma}
\begin{proof}
  In this proof, we refer to both a submatrix of
  $\mathcal{T}_j$ 
  and its determinant as a minor.  Let
  $M^{i_1,i_2,\ldots,i_r}_{j_1,j_2,\ldots,j_r}$ be a minor of 
  $\mathcal{T}_j$.  
  We first
  prove that this minor {\em is} trivially zero if and only if at
  least one of its diagonal entries is zero.  We then prove
  that this minor has at least one zero on its diagonal if and
  only if there is a $t\in \{1,\ldots ,r\}$ such that
  $$
  j_t > \Big\lceil\frac{i_t}{n-k}\Big\rceil k.
  $$
  
  To prove the first equivalence, suppose first that the minor
  $M^{i_1,i_2,\ldots,i_r}_{j_1,j_2,\ldots,j_r}$
  has at least one
  zero on its diagonal.  We denote the entry in row $m$ and
  column $n$ of $M^{i_1,i_2,\ldots,i_r}_{j_1,j_2,\ldots,j_r}$ by
  $(m,n)$.
  Note that, if $(m,n) = 0$, then $(m',n')
  = 0$ for all entries $(m',n')$ with $m'\leq m$ and $n'\geq n$.  Let
  $S_r$ denote the permutation group on r letters.  From the
  expression
\begin{equation}\label{eq1}
  \sum_{\sigma\in S_r}
  (\mbox{sgn}\, \sigma)(1,\sigma(1))(2,\sigma(2))\ldots(r,\sigma(r))
\end{equation}
giving the determinant of an $r\times r$ matrix over a
commutative ring, we see that this minor is zero whether viewed
as an element of $\F^{r\times r}$ or as an element of $R^{r\times
  r}$ in the manner defined above; in other words, it is
trivially zero.  We now prove the other direction.  It is easy to
see that every trivially zero $1\times 1$ and $2\times 2$ minor
of $\mathcal{T}_j$ has a zero on its diagonal.  Let $n$ be a
positive integer with $n\geq 3$, and suppose that, for $r\leq n -
1$, every trivially zero $r\times r$ minor has a zero on its
diagonal.  Suppose that the $n\times n$ minor
$M^{i_1,i_2,\ldots,i_n}_{j_1,j_2,\ldots,j_n}$ is trivially zero.
We use an induction argument to show that this minor has at least
one zero on its diagonal.  Notice that, when
$M^{i_1,i_2,\ldots,i_n}_{j_1,j_2,\ldots,j_n}$ is viewed as an
element of $R^{n\times n}$ in the manner defined above, the entry
$(n,1)$ appears exactly once.  If this entry is zero, then every
entry in the minor is zero, and thus all diagonal entries are
zero.  Suppose this entry is not zero.  Doing a cofactor
expansion along the first column shows that the $(n - 1)\times (n
- 1)$ minor $M^{i_1,i_2,\ldots,i_{n - 1}}_{j_2,\ldots,j_n}$ is
trivially zero.  By the induction hypothesis, this minor must
have a zero on its diagonal.  Thus, there is an entry $(s,s + 1)
= 0$ in $M^{i_1,i_2,\ldots,i_n}_{j_1,j_2,\ldots,j_n}$, which
means that the minor $M^{i_1,i_2,\ldots,i_s}_{j_{s +
    1},\ldots,j_n}$ is such that all of its entries are zero.
Because we assumed $M^{i_1,i_2,\ldots,i_n}_{j_1,j_2,\ldots,j_n}$
is trivially zero, it follows that at least one of the minors
$M^{i_1,i_2,\ldots,i_s}_{j_1,j_2,\ldots,j_s}$, $M^{i_{s +
    1},\ldots,i_n}_{j_{s + 1},\ldots,j_n}$ is trivially zero.  By
the induction hypothesis, at least one of these minors has a zero
on its diagonal.  As the diagonals of these minors lie on the
diagonal of $M^{i_1,i_2,\ldots,i_n}_{j_1,j_2,\ldots,j_n}$, we are
done.
  
  To prove the second equivalence, we simply note that the diagonal
  entries of this minor are the entries $(i_1,j_1)$, $(i_2,j_2)$,
  $\ldots$, $(i_r,j_r)$ of $\mathcal{T}_j$.  From the structure of
  $\mathcal{T}_j$, it is clear that the diagonal entry $(i_t,j_t)$
  is zero if and only if
$$
j_t > \Big\lceil\frac{i_t}{n-k}\Big\rceil k.
$$
\end{proof}

\begin{theo}                                \label{main1}
  Let $\C$ be an $(n,k,\delta)$ convolutional code described by
  matrices $(A,B,C,D)$ and consider the block Toeplitz matrix
  $\mathcal{T}_j$ introduced in\eqr{Bl-To}. Then $\C$ has $j$th
  column distance $d_j=(n-k)(j+1)+1$ if and only if every minor
  which is not trivially zero is nonzero.
\end{theo}

\begin{proof}
  $\Longleftarrow ~:$ Suppose that
  $$
  \left(
  \begin{array}{ccccccccc}
  y_0 &  y_1 &\ldots & y_j   &  \vline&u_0&u_1 & \ldots  &  u_j 
\end{array}\right)^T
$$
is a finite-weight codeword with $u_0\not = 0$ and that the
vector
$$
\left(
  \begin{array}{cccc}
  u_0 & u_1 &\ldots &  u_j 
  \end{array}
\right)^T
$$
has weight $r$.  Suppose that every minor of the
matrix\eqr{Bl-To} which is not trivially zero is nonzero.
Because $u_0\not = 0$, this means that at most $r-1$ rows
of\eqr{Bl-To} are in the left kernel of
$$
\left(
  \begin{array}{cccc}
  u_0 &u_1&  \ldots &  u_j 
  \end{array}
\right)^T.
$$
Thus, this codeword has weight at least $r + (j + 1)(n - k) -
(r - 1) = (j + 1)(n - k) + 1$.  We therefore have that the weight
of any codeword with $u_0\not = 0$ is at least $(j + 1)(n - k) +
1$.  In other words, $d_j \geq (j + 1)(n - k) + 1$.
Proposition~\ref{P-dcj.bound} implies that $d_j= (j + 1)(n - k) +
1$.
      
$\Longrightarrow ~:$ We prove the contrapositive.  We first note
that the result follows trivially if, for some integer 
$i$ with $0\leq i\leq j$, the matrix
$F_i$
contains a zero entry.  We
therefore assume that all such entries are nonzero.  Suppose
that the matrix\eqr{Bl-To} has an $r\times r$ minor which is zero
but not trivially zero, where $r\geq 2$ (in this proof, we again use 
'$r\times r$ minor' to refer to both the
submatrix and its determinant).  If the $r$ rows of this minor
belong to the left kernel of a column vector of weight $r$, then,
because of the structure of\eqr{Bl-To}, we can form a
nonzero vector
$$
\left(
  \begin{array}{cccc}
  u_0   & u_1&  \ldots&  u_j 
  \end{array}
\right)^T
$$
of weight $r$ with $u_0\not = 0$ such that the weight of
$$
\left(
  \begin{array}{ccccccccc}
  y_0 &  y_1 &  \ldots &  y_j  &   \vline&  u_0&
  u_1 &
  \ldots &   u_j
  \end{array}
\right)^T
$$
is at most $r + (j + 1)(n - k) - r = (j + 1)(n - k)$.  If not,
then the $r$ rows of this minor belong to the kernel of a nonzero
column vector of weight $r'\leq r$.  The $r'$ nonzero components
of this vector pick out an $r\times r'$ submatrix.  We would like
to see that this submatrix contains an $r'\times r'$ minor which
is zero but not trivially zero.  To obtain this minor, we simply
choose the bottom $r'$ rows.  This minor is clearly zero.
Because
we have assumed that the entries of the matrices $F_i$ are all
nonzero,
we know from
Lemma~\ref{B-minors2} that a minor is not trivially zero
if and only if its $i$th column has the property that the last $r
+ 1 - i$ entries are nonzero.  The columns of the original
$r\times r$ minor have this property.  Thus, the columns of our
$r'\times r'$ subminor have this property.  Because of the
structure of\eqr{Bl-To}, we can form a nonzero vector
$$
\left(
  \begin{array}{cccc}
  u_0   & u_1&  \ldots&  u_j 
  \end{array}
\right)^T
$$
of weight $r'$ with $u_0\not = 0$ such that the weight of
$$
\left(
  \begin{array}{ccccccccc}
  y_0 &  y_1 &  \ldots &  y_j  &   \vline&  u_0&
  u_1 &
  \ldots &   u_j
  \end{array}
\right)^T
$$
is at most $r' + (j + 1)(n - k) - r' = (j + 1)(n - k)$.
\end{proof}

Specializing Theorem~\ref{main1} for $j=L$, and
recalling~\cite[Corollary 2.4]{gl03r}, we immediately get an
algebraic criterion for an $(A,B,C,D)$ convolutional code to
represent a maximum distance profile code:

\begin{cor}                         \label{main-cor}
  Let $L=\Big\lfloor\frac{\delta}{k}
  \Big\rfloor+\Big\lfloor\frac{\delta}{n-k}\Big\rfloor$.  Then,
  the matrices $(A,B,C,D)$ generate a maximum distance profile
  $(n,k,\delta)$ convolutional code if and only if the matrix
  $\mathcal{T}_L$ has the property that every minor which is not
  trivially zero is nonzero.
\end{cor}

\begin{rem}
  Theorem~\ref{main1} and Corollary~\ref{main-cor} give
  polynomial conditions in the entries of the matrices
  $(A,B,C,D)$ which guarantee that a convolutional code has the
  maximum distance property. These algebraic properties are
  invariant under state space transformations. This means that,
  if $(A,B,C,D)$ has this property, then so does
  $(SAS^{-1},SB,CS^{-1},D)$ for every matrix $S\in
  Gl_\delta(\F)$. As a result, these conditions are really
  algebraic conditions on the quasi-projective
  variety~\cite{ha77p,ra94} describing the set of rate $k/n$
  convolutional codes of degree $\delta$. We have therefore
  established that the set of convolutional codes having the
  maximum distance property form a Zariski open set in this
  quasi-projective variety.
\end{rem}

The remainder of this section is devoted to showing that this
Zariski open set of maximal distance profile codes is nonempty as
soon as the base field is sufficiently large. As a first step
toward this result, we have the following theorem.
 
\begin{theo}                       \label{A}
  Let $j,k,n$ be fixed positive integers and consider the matrix
  introduced in\eqr{Bl-To}.  If the field $\F$ is sufficiently
  large, then one can find a sequence of matrices
  $\left\{F_0,\ldots,F_j\right\}$, where $F_i\in\F^{(n-k)\times
    k}~\forall i\in\left\{0,1,\ldots,j\right\}$, such that every
  minor of the matrix $ \mathcal{T}_j$ which is not trivially
  zero is nonzero.
\end{theo}
\begin{proof}
  Let $\F$ be an arbitrary finite field and let $\bar \F$ denote the
  algebraic closure of $\F$.  Note that $\bar \F$ is an infinite
  field.  To say that a minor of $\mathcal{T}_j$ is zero but not
  trivially zero is to say that the entries of $\mathcal{T}_j$
  satisfy a nonzero polynomial equation in $\bar\F
  [x_1,x_2,\ldots,x_{(j+1)(n-k)k}]$.  As there are finitely many
  minors, there are finitely many such polynomial equations
  describing those matrix sequences
  $\left\{F_0,\ldots,F_j\right\}$ for which $\mathcal{T}_j$ has
  at least one minor that vanishes but does not trivially vanish.
  Each of these polynomials describes a proper algebraic subset
  of $\bar\F^{(j+1)(n-k)k}$, the complement of which is a
  nonempty Zariski open set in $\bar\F^{(j+1)(n-k)k}$.  We may
  take the intersection of these Zariski open sets, and, as there
  are finitely many of them, the result is again a nonempty
  Zariski open set in $\bar\F^{(j+1)(n-k)k}$.  Take
  $\left\{F_0,\ldots,F_j\right\}$ to be an element in this
  intersection.  There are finitely many entries in this matrix
  sequence.  Thus, either all of the entries belong to $\F$, or
  there is a finite extension field of $\F$ containing all of
  them.  In this case, we may instead take $\F$ to be the
  smallest such extension.
\end{proof}

In order to show the existence of maximum distance profile
convolutional codes of arbitrary rate $k/n$, we use a result by
Tether~\cite{te70} from {\em minimal partial realization theory}.
For a given rate, this result allows us to show the existence of
such codes possessing only certain degrees.  It will require only
a small amount of additional work to obtain the existence of such
codes of arbitrary degree. Readers interested in the minimal
partial realization problem and its connections are referred
to~\cite{an86}.

\begin{theo}\label{B}
  Let $j,k$, and $n$ be fixed positive integers with $k<n$.  Let
  $\{F_0,F_1,\ldots,F_j\}$ and $\mathcal{T}_j$ be as in
  Theorem~\ref{A}.  Let $\delta$ be the smallest positive integer
  such that
  $$
  \Big\lfloor\frac{\delta}{k}\Big\rfloor +
  \Big\lfloor\frac{\delta}{n-k}\Big \rfloor\geq j.
  $$
  Suppose that $\{F_0,F_1,\ldots,F_j\}$ is such that the
  corresponding $\mathcal{T}_j$ has the property that its minors
  which are not trivially zero are nonzero.  Then, one can extend
  the sequence $\{F_0,F_1,\ldots,F_j\}$ to an infinite sequence
  $\{F_0,F_1,\ldots \}$ of $(n-k)\times k$ matrices over $\F$
  such that the infinite block Hankel matrix
  $$
  \mathcal{F}= \left(
     \begin{array}{cccc}
                     F_1    & F_2      & F_3     & \cdots  \\
                     F_2    & F_3      & F_4     & \cdots  \\
                     F_3    & \vdots   & \vdots  & \ddots  \\
                     \vdots & \vdots   & \vdots
     \end{array}
   \right)
   $$
   has rank $\delta$.
\end{theo}
\begin{proof}
  Let $\mathcal{F}_{x,y}$ denote the block Hankel matrix
  $$
  \mathcal{F}_{x,y}= \left(
     \begin{array}{cccc}
                     F_1    & F_2    & \cdots & F_y      \\
                     F_2    & F_3    & \cdots & F_{y+1}  \\
                     \vdots & \vdots & \cdots & \vdots   \\
                     F_x    & F_{x+1}& \cdots & F_{x+y-1}
     \end{array}
   \right)
   $$
   As was shown in~\cite[Theorem 1]{te70}, any matrix sequence
   $\{F_1,\ldots,F_j\}$ has a minimal partial realization of
   degree $d$, where
   $$
   d = \sum_{i=1}^j\rank \mathcal{F}_{i,j+1-i}
   -\sum_{i=1}^{j - 1}\rank \mathcal{F}_{i,j-i}.
   $$
   It is easy to see that, up to a reordering of block
   columns, each $\mathcal{F}_{i,j+1-i}$ appearing in the first
   summation in the above formula for $d$ is a submatrix of
   $\mathcal{T}_j$ (take the intersection of the last $i$ block
   rows and the first $j+1-i$ block columns).  As each
   $\mathcal{F}_{i,j-i}$ in the second summation is a submatrix
   of $\mathcal{F}_{i,j+1-i}$, the same is true for these.  By
   assumption, the only minors of $\mathcal{T}_j$ that are zero
   are those that are trivially zero.  Thus, $\rank
   \mathcal{F}_{x,y} = \min ((n-k)x, ky)$ for each
   $\mathcal{F}_{x,y}$ appearing in the above formula for $d$.
   
   First, suppose that $k\geq n-k$.  The formula for $d$ then
   becomes
   $$
   d=\sum_{i=1}^j\min (i(n-k),(j+1-i)k)-\sum_{i=1}^{j-1}\min
   (i(n-k),(j-i)k).
   $$
   Suppose there exists an integer $r, 1\leq r\leq j-1$, with
   $r(n-k)\geq (j-r)k$.  Let $i^*$ be the smallest such integer.
   Then, the formula for $d$ becomes
   $$
   d=(i^*-1)(n-k) +(\min (i^*(n-k),(j-i^*+1)k)-\min
   ((i^*-1)(n-k),(j-i^*+1)k)).
   $$
   By definition of $i^*$, $(i^*-1)(n-k)<(j-i^*+1)k$, so that
   the last term in this expression is $-(i^*-1)(n-k)$.  Thus, 
$$
   d = \min (i^*(n-k),(j-i^*+1)k).
$$  
   Note that
   $i^*=\lceil j \frac{k}{n} \rceil$.  Consider the difference
\begin{multline*}
  (j-i^*+1)k-(i^*-1)(n-k)= (j-\Big\lceil j\frac{k}{n} \Big\rceil
  +1)k-
  (\Big\lceil j\frac{k}{n} \Big\rceil -1)(n-k)=\\
  (j+1)k-\Big\lceil j\frac{k}{n} \Big\rceil n +(n-k).
\end{multline*}
We want to see that this number is at least $n-k$, or,
equivalently, that $(j+1)k-\lceil j\frac{k}{n} \rceil n\geq 0$.
This will imply that $\min (i^*(n-k),(j-i^*+1)k)=i^*(n-k)$.  We
have

$$
(j+1)k-\Big\lceil j\frac{k}{n} \Big\rceil n\geq 0
\Longleftrightarrow (j+1)k\geq \Big\lceil j\frac{k}{n} \Big\rceil
n \Longleftrightarrow \frac{(j+1)k}{n}\geq \Big\lceil
j\frac{k}{n} \Big\rceil.
$$

By assumption, $k\geq n-k$, which means $k\geq \frac{n}{2}$.  It
follows that $\frac{(j+1)k}{n}\geq \lceil j\frac{k}{n} \rceil$.
Thus, we have $d=i^*(n-k)$.

By assumption, $i^*(n-k)\geq (j-i^*)k$.  We have
\begin{multline*}
  i^*(n-k)\geq (j-i^*)k \Longleftrightarrow i^*(n-k)+i^*k\geq jk
  \Longleftrightarrow \frac{i^*(n-k)}{k}+i^*\geq j\\
  \Longleftrightarrow \Big\lfloor
  \frac{i^*(n-k)}{k}+i^*\Big\rfloor =\Big\lfloor
  \frac{i^*(n-k)}{k}\Big\rfloor +\Big\lfloor i^*\Big\rfloor
  =\Big\lfloor\frac{i^*(n-k)}{k} \Big\rfloor
  +\Big\lfloor\frac{i^*(n-k)}{n-k} \Big\rfloor\geq j.
\end{multline*}
Let $\delta=i^*(n-k)$.  We want to see that $\delta$ is the
smallest positive integer satisfying
$\lfloor\frac{\delta}{k}\rfloor
+\lfloor\frac{\delta}{n-k}\rfloor\geq j$.  Consider
$\lfloor\frac{\delta}{k}\rfloor
=\lfloor\frac{i^*(n-k)}{k}\rfloor$.  Since $i^*=\lceil
j\frac{k}{n}\rceil$, we may write $i^*=\frac{jk+s}{n}$ where
$s\in\N_0$ and $0\leq s\leq n-1$.  Then,
$\lfloor\frac{\delta}{k}\rfloor = \frac{znk-y}{nk}$ where
$z,y\in\N_0$ and $0\leq y\leq nk-1$.  Thus, we have
$$
\Big\lfloor\frac{\delta}{k}\Big\rfloor
+\Big\lfloor\frac{\delta}{n-k}\Big\rfloor =\frac{znk-y}{nk}
+\frac{jk+s}{n}\leq \frac{(jk+s)(n-k)}{nk} +\frac{jk+s}{n}
$$
with equality precisely when $y=0$ or, in other words, when
$k~\vline ~i^*(n-k)$.  Suppose $y=0$.  Then, since $s\leq n-1$
and $k\geq \frac{n}{2}$, we have
$$
\Big\lfloor\frac{\delta}{k}\Big\rfloor
+\Big\lfloor\frac{\delta}{n-k}\Big\rfloor =\frac{(jk+s)(n-k)}{nk}
+\frac{jk+s} {n} =\frac{jk+s}{k} =j+\frac{s}{k}<j+2
$$
As $j+\frac{s}{k}$ must be an integer, we may write
$j+\frac{s}{k}\leq j+1$.  For the same reason, we must have
$s=lk, l\in \{0,1\}$.  From this, we see that
$\lfloor\frac{\delta}{k}\rfloor
=\lfloor\frac{i^*(n-k)}{k}\rfloor$ is bounded above by $j+1$, and
this upper bound is obtained precisely when $k~\vline ~i^*(n-k)$
and $l=1$.  Because $(n-k)~\vline ~\delta$, it follows that when
$l=1$, $\lfloor\frac{\delta -1}{k}\rfloor +\lfloor\frac{\delta
  -1}{n-k}\rfloor\leq j-1$, and when $l=0$, that
$\lfloor\frac{\delta -1}{k}\rfloor +\lfloor\frac{\delta
  -1}{n-k}\rfloor\leq j-2$.  If $y$ is nonzero, then 
$\lfloor\frac{\delta}{k}\rfloor
+\lfloor\frac{\delta}{n-k}\rfloor\leq j$, and so 
$\lfloor\frac{\delta -1}{k}\rfloor
+\lfloor\frac{\delta -1}{n-k}\rfloor\leq j-1$.

If no such $r$ exists, there are two possibilities.  The first is
that $j(n-k)>k$, so that $\min (j(n-k),k)=k$.  Then
$d=(j-1)(n-k)+k-(j-1)(n-k)=k$.  Because $(j - 1)(n - k) < k$, we
have the inequalities
$\frac{n}{n-k} =\frac{k}{n-k} +1>j>\frac{k}{n-k}$, and so we may
write $\lfloor\frac{k}{n-k}\rfloor +1\geq
j>\lfloor\frac{k}{n-k}\rfloor$.  Letting $\delta =k$, we see that
$\delta$ is the smallest positive integer satisfying
$\lfloor\frac{\delta}{k}\rfloor
+\lfloor\frac{\delta}{n-k}\rfloor\geq j$.  The second possibility
is that $j(n-k)\leq k$, so that $\min (j(n-k),k)=j(n-k)$.  Let
$\delta =j(n-k)$.  Clearly, $\delta$ is the smallest positive
integer satisfying $\lfloor\frac{\delta}{k}\rfloor
+\lfloor\frac{\delta}{n-k}\rfloor\geq j$.
  
The proof for the case $k<n - k$ is similar.
\end{proof}

Theorem~\ref{B} immediately establishes the existence of maximum
distance profile codes for certain parameters $(n,k,\delta)$:
\begin{lemma}\label{C}
  Let $k,n$ and $\delta$ be positive integers such that $k<n$ and
  either $k \, | \,\delta$ or $n-k\, | \,\delta$.  Then, an
  $(n,k,\delta)$ maximum distance profile convolutional code
  exists over a sufficiently large base field.
\end{lemma}
\begin{proof}
  Let $\lfloor\frac{\delta}{k}\rfloor +
  \lfloor\frac{\delta}{n-k}\rfloor = L$.  Note that if $\delta$
  decreases, then $L$ must decrease.  From the proof of
  Theorem~\ref{A}, we know that there exists a finite extension
  field of $\F$ and a sequence $\left\{F_0,
    F_1,\ldots,F_L\right\}$ of $(n-k)\times k$ matrices with
  entries in this extension field such that the only minors of
  the corresponding $\mathcal{T}_L$ which are zero are those
  which are trivially zero.  From the proof of Theorem~\ref{B},
  we know that there exist matrices $A \in \F^{\delta \times
    \delta}, \,\,\, B \in \F^{\delta \times k}$, and $C \in
  \F^{(n-k) \times \delta}$ which give a minimal partial
  realization of the sequence $\left\{F_1,\ldots,F_L\right\}$.
  Let $D = F_0$.  Then, the matrices $A,B,C$, and $D$ describe an
  $(n,k,\delta)$ convolutional code via\eqr{iso}.  By
  Corollary~\ref{main-cor}, this code has a maximum distance
  profile.
\end{proof}

We now state and prove the main Theorem.
\clearpage

\begin{theo}             \label{main}
  Let $k$ and $n$ be positive integers such that $k<n$.  Let
  $\delta$ be a positive integer.  Then, there exists a maximum
  distance profile $(n,k,\delta)$ convolutional code over some
  finite extension field of $\F$.  Moreover, the set of maximum
  distance profile $(n,k,\delta)$ convolutional codes forms a
  generic set in $\bar\F^{(\delta + n - k)(\delta - k)}$.
\end{theo}
\begin{proof}
  Let $L = \lfloor\frac{\delta}{k}\rfloor +
  \lfloor\frac{\delta}{n-k}\rfloor$.  Let $\delta^*$ be the
  smallest integer satisfying this equality.  If $L = 0$, the
  theorem is easily seen to be true.  This is because the problem
  of proving existence is reduced to finding an $(A,B,C,D)$
  representation of an $(n,k,\delta)$ convolutional code over a
  finite extension field of $\F$ such that $[-I \ D]$ represents
  the parity check matrix of an MDS block code. The proof that
  such a $D$ can be found is essentially identical to the proof
  of Lemma~\ref{A}.  The set of such $(A,B,C,D)$ representations
  is obviously a Zariski open set in $\bar\F^{(\delta + n -
    k)(\delta - k)}$.  This proves that the set of such codes
  forms a generic set in $\bar\F^{(\delta + n - k)(\delta - k)}$.
  Thus, we may assume $L\not = 0$.
  
  Let $S$ denote the set of 4-tuples of matrices $(A,B,C,D)$
  with $A \in \F^{\delta\times \delta}, \,\,\, B \in 
  \F^{\delta \times k}, \,\,\, C \in\F^{(n-k) \times \delta}$,
  and 
  $D \in \F^{(n-k) \times k}$ and having the property that every
  minor in\eqr{Bl-To} which is not trivially zero is nonzero.
  We first show $S$ is a
  nonempty Zariski open set in $\bar\F^{(\delta + n - k)
  (\delta - k)}$.  $S$ is obviously Zariski open.  
  We now show it is nonempty.  Let $r = \delta - \delta^*$.  
  Let $(\tilde A,\tilde B,\tilde C,\tilde D)$ be a representation 
  of a maximum distance profile $(n,k,\delta ^*)$ convolutional
  code; the existence of such a code is implied by
  Lemma~\ref{C}.  Consider the matrices
$$
  A = \left(\begin{array}{c|c}
      0_{r\times r}     & 0_{r\times \delta^*} \\ 
      \hline \\
      0_{\delta^*\times r}    & \tilde A                    
             \end{array}\right), \ 
         B = \left(\begin{array}{c}
        0_{r\times k}      \\
          \hline                 \\ 
               \tilde B                   
            \end{array}\right), \ 
           C = \left(\begin{array}{c|c}
               0_{(n - k)\times r}                     
              \end{array} 
             \begin{array}{r}
             \tilde C\end{array}\right), \ 
  D = \tilde D.
$$
  Notice that $CA^{i - 1}B = \tilde C\tilde A^{i - 1}\tilde B$
  for all $i\geq 1$.  Because $(\tilde A,\tilde B,\tilde C,\tilde
  D)$ is a representation of a maximum distance profile
  convolutional code, we have shown that $S$ is 
  nonempty.
  
  To complete the proof, we note that the reasoning used in the
  proof of Lemma~\ref{A} implies that 4-tuples of matrices
  $(A,B,C,D)$ in $\bar \F^{(\delta + n - k)(\delta - k)}$ such
  that $(A ,B)$ is a controllable pair and $(A ,C)$ is an
  observable pair form a nonempty Zariski open set in $\bar
  \F^{(\delta + n - k)(\delta - k)}$.  Intersecting this set with
  $S$ gives a nonempty Zariski open set in $\bar \F^{(\delta + n
    - k)(\delta - k)}$ consisting of maximum distance profile
  $(n,k,\delta)$ convolutional codes.  Thus, these codes form a
  generic set in $\bar \F^{(\delta + n - k)(\delta - k)}$.
\end{proof}

It was pointed out above that, when $n - k \, | \, \delta$, an
$(n,k,\delta)$ convolutional code has a maximum distance profile
if and only if it is strongly MDS.  With this, we have the
following Corollary:
\begin{cor}
When $n - k \, | \, \delta$, there exists an $(n,k,\delta)$
strongly MDS convolutional code over a finite extension field of $\F$.
\end{cor}

\section{Codes with Maximum Distance Profile in Terms of 
  Polynomial Generator Matrices}

In the coding literature, convolutional codes are usually studied
via (polynomial) generator and parity check matrices.  The
relevant results presented in~\cite{gl03r} were formulated in
terms of such polynomial matrix descriptions.  In this section,
we make the connection between polynomial and state space
descriptions of convolutional codes.  For this, we
follow~\cite{ro96a1,ro99a}, where further details may be found as
well.  We also state a necessary and sufficient condition for a
polynomial parity check matrix to define a maximum distance
profile convolutional code.

Consider the transfer function $T(z) := C(zI - A)^{-1} B + D$.
Let $P(z)^{-1} Q(z) = T(z)$ be a left coprime factorization of
$T(z)$ and $H(z) := [P(z)~ Q(z)]$.  Consider the polynomial
vectors:
$$
u(z) = u_0 z^{\gamma} + u_1 z^{\gamma - 1} + \hdots +
u_{\gamma};~ u_t\in \F^k , t = 0,\ldots ,\gamma,
$$
and
$$
y(z) = y_0 z^{\gamma} + y_1 z^{\gamma - 1} + \hdots +
y_{\gamma};~ y_t\in \F^{n-k} , t = 0,\ldots , \gamma.
$$
Then, the following conditions are equivalent.
\begin{enumerate}
\item The vectors $u_t$ and $y_t$ satisfy the state space
  equation\eqr{iso}.
\item The vectors $u_t$ and $y_t$ satisfy
\begin{equation}                   
\left(
\begin{array}{ccc|ccccc}
&  & &      D        &               &      &  &  \\
&  & &     CB        &        D      &      &  &  \\
&-I&\ \ &     CAB    &        CB     &\ddots&  &  \\
&  & &        \vdots &               &\ddots&\ddots &  \\
&  & &CA^{\gamma -1}B&CA^{\gamma -2}B&\cdots& CB & D
\end{array}
\right) \left(
\begin{array}{c}
y_0 \\
y_1 \\
\vdots  \\
y_\gamma  \\   \hline
u_0\\
u_1 \\
\vdots  \\
u_\gamma
\end{array}
\right) =0.  
\end{equation}
\item There exists a `state vector'
  $$
  x(z) = x_0z^{\gamma}+x_{1}z^{{\gamma}-1} + \ldots +
  x_{\gamma};\ x_t\in\F^\delta, t=0,\ldots,\gamma,
  $$
  such that
\begin{equation}
  \label{kern}
\left[
   \begin{array}{ccc}
   zI-A&0_{\delta\times (n-k)}&-B\\ -C&I_{n-k}&-D
   \end{array}
 \right]\left[
   \begin{array}{c}
   x(z)\\y(z)\\u(z)
   \end{array}
 \right] =0.
\end{equation}
\item $\zwei{y(z)}{u(z)}$ is a code word, i.e.
  $$
  H(z)\zwei{y(z)}{u(z)}= [P(z) \ Q(z)] \zwei{y(z)}{u(z)}=0.
  $$
\item $y(z)=T(z)u(z)$
\end{enumerate}
The proof of these equivalences is straightforward, and more
details can be found, e.g., in~\cite{ro99a}.

We close with the following theorem.
\begin{theo}
  Let $H(z) = \sum_{l = 0}^{\mu} H_l z^l$ be the parity check
  matrix of an $(n,k,\delta)$-code.\\ Assume $H_l = 0$ for $l >
  \mu$.  Let
  $$
  \mathcal{H}_j := \left(\begin{array}{cccc}
      H_0                                                 \\
      H_1     & H_0                                       \\
      \vdots  & \vdots    & \ddots                        \\
      H_j & H_{j-1} & \cdots & H_0\end{array} \right) \in \F^{(j
    + 1)(n - k)\times (j + 1)n},
  $$
  Then $H(z)$ represents a code whose $jth$ column distance
  $d_j = (n - k)(j + 1)+ 1$ if and only if every $(j + 1)(n -
  k)\times (j + 1)(n - k)$ full-size minor formed from the
  columns with indices $1\leq i_1 <\ldots < i_{(j + 1)(n - k)}$,
  where $i_{s(n - k)}\leq sn$ for $s = 1,\ldots , j$, is nonzero.
  
  In particular when
  $$
  j := L := \Big\lfloor\frac{\delta}{k}\Big\rfloor +
  \Big\lfloor\frac {\delta}{n - k}\Big\rfloor.
  $$
  then $H(z)$ represents a maximum distance profile code if
  and only if every $(L + 1)(n - k)\times (L + 1)(n - k)$
  full-size minor formed from the columns with indices $1\leq i_1
  <\ldots < i_{(L + 1)(n - k)}$, where $i_{s(n - k)}\leq sn$ for
  $s = 1,\ldots , L$, is nonzero.
\end{theo}
\begin{proof}
  This theorem is a direct consequence of~\cite[Corollary
  2.4]{gl03r} and ~\cite[Theorem 5.3]{gl03r}.
\end{proof}

\Section{Conclusion}
In this paper, we established the existence of maximum distance
profile codes for all transmission rates and all degrees. The
main results are existence results. Important questions remain
open as to how  maximum distance profile codes may be constructed
and the minimal field size required for doing so. For the
construction of such codes we found that many cyclic
convolutional codes~\cite{gl02u}  have the maximum distance
profile property and this might be a promising avenue for
constructing such codes.

The properties of these codes are very appealing for error
control coding; the distance between two trajectories which start
at a common initial state is maximal, and hence these codes have
the potential to have the maximal amount of errors per time
interval corrected. In applications where fault-diagnosis is
important (see e.g.~\cite{fl02,ha03}), it has already been
pointed out that codes with maximal free distance and hence also
codes with a maximal distance profile are very important.\bigskip

\noindent
{\bf Acknowledgement:} The authors wish to thank Heide
Gluesing-Luerssen for helpful comments throughout the preparation
of this paper.


\begin{thebibliography}{10}

\bibitem{an86}
A.C. Antoulas.
\newblock On recursiveness and related topics in linear systems.
\newblock {\em IEEE Trans. Automat. Contr.}, AC-31(12):1121--1135, 1986.

\bibitem{fl02}
M.~Fliess.
\newblock On the structure of linear recurrent error-control codes.
\newblock {\em ESAIM Control Optim. Calc. Var.}, 8:703--713 (electronic), 2002.
\newblock A tribute to J. L. Lions.

\bibitem{gl03r}
H.~Gluesing-Luerssen, J.~Rosenthal, and R.~Smarandache.
\newblock Strongly {MDS} convolutional codes, March 2003.
\newblock E-print math.RA/0303254.

\bibitem{gl02u}
H.~Gluesing-Luerssen and W.~Schmale.
\newblock On cyclic convolutional codes, October 2002.
\newblock E-print math.RA/0211040.

\bibitem{ha03}
C.~N. Hadjicostis.
\newblock Non-concurrent error detection and correction in fault-tolerant
  linear finite-state machines.
\newblock {\em IEEE Trans. Automat. Contr.}, 2003.
\newblock To appear.

\bibitem{ha77p}
M.~Hazewinkel.
\newblock Moduli and canonical forms for linear dynamical systems {III}: The
  algebraic geometric case.
\newblock In {\em Proc. of the 76 Ames Research Center ({NASA}) Conference on
  Geometric Control Theory}, pages 291--336. Math.Sci. Press, 1977.

\bibitem{jo89}
R.~Johannesson and K.~Zigangirov.
\newblock Distances and distance bounds for convolutional 
 codes -- an overview.
\newblock In {\em Topics in Coding Theory. In honour of L. H. Zetterberg.},
  Lecture Notes in Control and Information Sciences 
  \# 128, pages 109--136.
  Springer Verlag, 1989.

\bibitem{jo99}
R.~Johannesson and K.~Sh. Zigangirov.
\newblock {\em Fundamentals of Convolutional Coding}.
\newblock IEEE Press, New York, 1999.

\bibitem{ra94}
M.~S. Ravi and J.~Rosenthal.
\newblock A smooth compactification of the space of transfer functions with
  fixed {McM}illan degree.
\newblock {\em Acta Appl. Math}, 34:329--352, 1994.

\bibitem{ro96a1}
J.~Rosenthal, J.~M. Schumacher, and E.~V. York.
\newblock On behaviors and convolutional codes.
\newblock {\em IEEE Trans. Inform. Theory}, 42(6, part 1):1881--1891, 1996.

\bibitem{ro99a1}
J.~Rosenthal and R.~Smarandache.
\newblock Maximum distance separable convolutional codes.
\newblock {\em Appl. Algebra Engrg. Comm. Comput.}, 10(1):15--32, 1999.

\bibitem{ro99a}
J.~Rosenthal and E.~V. York.
\newblock {BCH} convolutional codes.
\newblock {\em IEEE Trans. Inform. Theory}, 45(6):1833--1844, 1999.

\bibitem{sm01a}
R.~Smarandache, H.~Gluesing-Luerssen, and J.~Rosenthal.
\newblock Constructions for {MDS}-convolutional codes.
\newblock {\em IEEE Trans. Inform. Theory}, 47(5):2045--2049, 2001.

\bibitem{te70}
A.~J. Tether.
\newblock Construction of minimal linear state--variable models from finite
  input--output data.
\newblock {\em IEEE Transactions on Automatic Control}, 15:427--436, 1970.

\end{thebibliography}

\end{document}